\documentclass[twocolumn]{svjour3}
  \def\makeheadbox{\relax}
 
\smartqed  %
\usepackage{graphicx}
\usepackage{amsmath, amssymb}
\usepackage{latexsym}
\usepackage{mathtools}

\providecommand{\Div}{\operatorname{div}}          %
\providecommand{\Tr}{\operatorname{tr}}                    %

\newcommand{\Vg}{{\mathbf{g}}}

\newcommand{\Vn}{{\mathbf{n}}}

\newcommand{\Vu}{{\mathbf{u}}}

\newcommand{\Vx}{{\mathbf{x}}}

\providecommand{\Bx}{{\boldsymbol{x}}}

\newcommand{\VD}{{\mathbf{D}}}

\newcommand{\VF}{{\mathbf{F}}}

\newcommand{\VP}{{\mathbf{P}}}

\newcommand{\VT}{{\mathbf{T}}}

\newcommand{\VV}{{\mathbf{V}}}

\newcommand{\Jmat}{{\rm J}}

\providecommand{\Ci}{{\cal I}}

\providecommand{\Cm}{{\cal M}}

\providecommand{\Cu}{{\cal U}}

\providecommand{\bbN}{\mathbb{N}}

\providecommand{\bbR}{\mathbb{R}}

\newcommand{\DX}{\,\mathrm{d}\Bx}

\newcommand{\J}{\Jmat}
\newcommand{\dJ}{\mathrm{d\J}}
\newcommand{\ddJ}{\mathrm{d^2\J}}
\newcommand{\T}{\VT}
\newcommand{\F}{\VF_i}

\newcommand{\dT}{\VV}%

\newcommand{\ders}[1]{\mathrm{d_s}(#1)}

\newcommand{\myref}[1]{\hat#1}%
\newcommand{\Kref}{\myref{K}}
\newcommand{\xref}{{\myref{\Bx}}}
\newcommand{\DXref}{\,\mathrm{d}\xref}
\usepackage{tikz}
\usetikzlibrary{shapes,arrows,positioning}
\usepackage{cite}
\usepackage{newfloat}
\DeclareFloatingEnvironment[name=Lst.]{mylisting}
\usepackage{hyperref}

\begin{document}

\title{Automated shape differentiation in the Unified Form Language}
\author{David A. Ham \and Lawrence Mitchell \and Alberto Paganini \and Florian Wechsung}
\authorrunning{D. A. Ham, L. Mitchell, A. Paganini, and F. Wechsung} %
\institute{\and D. A. Ham; Department of Mathematics, Imperial College London, London, SW7 2AZ, UK;\\
             \email{david.ham@imperial.ac.uk}
             \\
             \and
             Lawrence Mitchell; Department of Computer Science, Durham University, Durham, DH1 3LE, UK.\\
             \email{lawrence.mitchell@durham.ac.uk}
             \\
                         \and 
             A. Paganini and F. Wechsung \at
             Mathematical Institute, University of Oxford,\\ Oxford, OX2 6GG, UK. \\
             \email{paganini@maths.ox.ac.uk, wechsung@maths.ox.ac.uk}
}

\date{}

\maketitle

\begin{abstract}
We discuss automating the {calculation} of weak shape derivatives in
the Unified Form Language (Aln{\ae}s et al., ACM
Trans. Math. Softw., 2014) by introducing an
appropriate additional step in the pullback from physical to
reference space that computes G\^{a}teaux derivatives with respect
to the coordinate field.  We illustrate the ease of use with several
examples.

\end{abstract}

\section{Introduction}

Physical models often involve functionals that assign real values to the
solutions of partial differential equations (PDEs).
For instance, the compliance of a structure is a function of the solution
to the elasticity equations, and the drag of
a rigid obstacle immersed in a fluid is a function of the
solution to the Navier-Stokes equations. 

This type of functional depends on the PDE parameters, and it is often
possible to compute the derivative of a functional with respect to a
chosen set of parameters.  This derivative can in turn be used to
perform sensitivity analysis and to run optimization algorithms with
respect to parameters in the PDE.

The shape of the physical domain that is part the PDE-model
(like the shape of the rigid obstacle mentioned above) is a parameter
that is not straightforward to handle. Although we can compute the
shape derivative of a functional following shape calculus rules \cite{DeZo11},
this differentiation exercise is often tedious and error prone.
In \cite{Sc18}, Schmidt introduces the open-source library FEMorph:
an automatic shape differentiation toolbox for the Unified Form Language
(UFL, \cite{AlLoOlRoWe14}). \textcolor{black}{The library }FEMorph is based on refactoring UFL expressions
and applying shape calculus differentiation rules recursively.
It can compute first- and second-order shape derivatives
(both in so-called weak and strong form), and it has been successfully employed
to solve shape optimization problems \cite{ScScWa18}.

This \textcolor{black}{article} presents an alternative approach to automated
shape differentiation. The key idea is to rely solely on pullbacks
and standard G\^ateaux derivatives. This approach is more generic and robust,
because it does not require handling of special cases.
In particular, it circumvents the implementation of shape
calculus rules, and required only a minor modification of UFL, because
UFL supports G\^ateaux derivatives with respect to functions.
As a result, UFL is now capable of shape differentiating any integral that can be
be expressed in it. 

\textcolor{black}{This article is organized as follows.
In Section \ref{sec:theory},
we revisit shape calculus and describe how to shape differentiate using
standard finite element software. In Section \ref{sec:examples},
we consider three test cases and show how to compute shape derivatives using
Firedrake and UFL. In Section \ref{sec:implementation-validation},
we describe code validation of this new UFL feature.
In Section \ref{sec:shapeoptexample}, we solve a PDE-constrained shape optimization
test case with a descent algorithm implemented in Firedrake and UFL.
Finally, in Section \ref{sec:discussion},
we summarize the contribution of this article.}

\section{Shape differentiation on the reference element}
\label{sec:theory}

A \emph{shape functional} is a map
$\J:{\Cu}\to\bbR$ defined on the collection of domains
${\Cu}$ in $\bbR^d$, where $d\in\bbN$ denotes the space dimension.
{We follow the perturbation of identity approach \cite{simon_differentiation_1980} and for} %
a vector field  $\VV\in W^{1,\infty}({\bbR^d},\bbR^d)$,
we consider the family of transformations $\{\VP_s\}_{s\geq0}$
defined by {$\VP_s(\Vx) = \Vx + s\dT(\Vx)$ for every $\Vx\in\bbR^d$}. Note that
$\VP_s$ is a diffeomorphism for any sufficiently small $s\in\bbR^+$.

{For a domain $\Omega\in {\Cu}$}, let $\Omega_s\coloneqq \VP_s(\Omega)$ and assume that $\Omega_s\in {\Cu}$ for $s$ sufficiently small.
The \emph{shape directional derivative}
of $\J$ at $\Omega$ in direction $\VV$ is the derivative
\begin{equation*}
\dJ(\Omega)[\VV]\coloneqq \lim_{s\searrow 0}\frac{\J(\Omega_s)-\J(\Omega)}{s}\,.
\end{equation*}

We say that  $\J$ is \emph{shape differentiable in $\Omega$} if the directional derivative $\dJ(\Omega)[\VV]$ exists for every direction
$\VV\in W^{1,\infty}({\bbR^d},\bbR^d)$, and if
the \textcolor{black}{associated} map $\dJ(\Omega): W^{1,\infty}({\bbR^d},\bbR^d)\mapsto \bbR$ is linear and continuous.
In this case, the linear operator $\dJ(\Omega)$ is called the \emph{shape derivative} of $\J$ in $\Omega$.

To illustrate the shape differentiation of a shape functional, we consider
the prototypical example
\begin{equation}\label{eq:proto-j}
    \J(\Omega_s) = \int_{\Omega_s} u_{\VP_s} \DX\,,
\end{equation}
where $u_{\VP_s}$ is a scalar function\footnote{Shape functionals that involve vector fields
or boundary integrals can be treated following the same steps and employing
suitable pullbacks.}.
The subscript $\VP_s$ highlights the possible dependence of $u_{\VP_s}$ on the domain $\Omega_s$.

The standard procedure to compute $\dJ$ is to employ
transformation techniques and rewrite
\begin{equation}\label{eq:JPs}
    \J(\Omega_s)= \int_{\Omega} (u_{\VP_s}\circ\VP_s) \det(\VD\VP_s)\DX\,,
\end{equation}
{where  $\VD\VP_s$ denotes the Jacobian matrix of $\VP_s$,
and $u_{\VP_s}\circ\VP_s$ denotes the composition of $u_{\VP_s}$
with $\VP_s$, that is, $(u_{\VP_s}\circ\VP_s)(\Vx) =u_{\VP_s}(\VP_s(\Vx))$ for every $\Vx \in \Omega$.}
Note that $\det(\VD\VP_s)>0$ for $s$ sufficiently small.
Then, by linearity of the integral,
the shape derivative $\dJ$ is given by
\begin{align}
\nonumber
    \dJ(\Omega)[\VV] &= \int_{\Omega}  \ders{(u_{\VP_s}\circ\VP_s)\det\VD\VP_s}\DX \\
\label{eq:dJgen}
    &= \int_\Omega \ders{u_{\VP_s}\circ \VP_s} + u_{\VP_0} \Div(\VV)\DX\,,
\end{align}
where $\ders{{}\cdot{}}$ denotes the derivative with respect to $s$ at $s=0$.
The term $\ders{u_{\VP_s}\circ \VP_s}$ is often called the \emph{material derivative} \cite{Be10}.
Its explicit formula depends on whether
the function $u_{\VP_s}$ does or does not dependent on $\Omega_s$ (see Section \ref{sec:examples}).

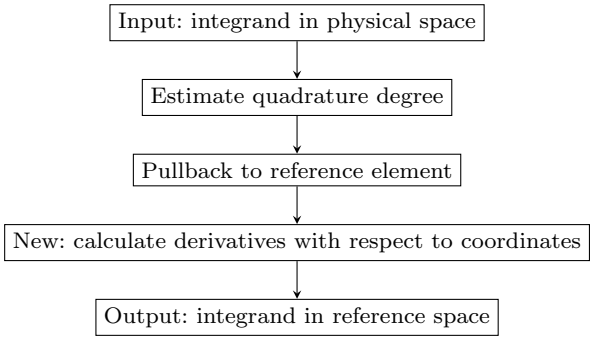
\begin{figure}[htbp]
    \begin{center}
        \tikzstyle{block} = [rectangle, draw,
        text centered, minimum height=1em,
        node distance=0.5cm]
        \begin{tikzpicture}
            \node[block] (function-deriv) {Input: integrand in physical space};
            \node[block, below=of function-deriv] (degree) {Estimate quadrature degree};
            \node[block, below=of degree] (pullback) {Pullback to reference element};
            \node[block, below=of pullback] (coord-deriv) {New: calculate derivatives with respect to coordinates};
            \node[block, below=of coord-deriv] (assemble) {Output:
              integrand in reference space};
            \draw[-stealth] (function-deriv) -- (degree);
            \draw[-stealth] (degree) -- (pullback);
            \draw[-stealth] (pullback) -- (coord-deriv);
            \draw[-stealth] (coord-deriv) -- (assemble);
        \end{tikzpicture}
        \caption{Symbolic workflow in UFL to transform integrals from
          physical to reference space.}
          \label{fig:UFLworkflow}
    \end{center}
\end{figure}

Next, we repeat the derivation of \eqref{eq:dJgen} in the context of finite elements
and derive an alternative formula for $\dJ$.
Let $\{K_i\}_{i\in \Ci}$ be a partition of $\Omega$
such that $\dot\cup_i \overline K_i = \overline{\Omega}$ and such that
the elements $K_i$ are non-overlapping. Additionally, let
$\{\F\}_{i\in \Ci}$ be a family of diffeomorphisms such that
$\F(\Kref) = K_i$ for every $i\in\Ci$,
where $\Kref$ denotes a reference element.
This induces a partition $\{\VP_s(K_i)\}_{i\in \Ci}$ of $\Omega_s$.
To evaluate \eqref{eq:proto-j}, standard finite element software
rewrites it as
\begin{align}
\nonumber
\J(\Omega_s) 
&= \sum_{i \in \Ci}\int_{\VP_s(K_i)}   u_{\VP_s} \DX \\
\label{eq:JOsFEM}
&= \sum_{i \in \Ci}\int_{\Kref}  (u_{\VP_s} \circ \VP_s\circ \F)  \vert\det(\VD (\VP_s\circ \F))\vert\DXref\,.
\end{align}
{Let $\VF_i^{-1}$ denote the inverse of $\VF_i$, that is, $\VF_i^{-1}(\VF_i(\Vx))= \Vx$
for every $\Vx\in\Kref$, and $\VF_i(\VF_i^{-1}(\Vx))= \Vx$ for every $\Vx\in K_i$.}
Since $\VP_s = (\VP_s\circ\VF_i)\circ\VF_i^{-1}$ and $\VP_s\circ\VF_i = \VF_i + s\VV\circ\VF_i$,
we can rewrite equation \eqref{eq:JOsFEM} as follows
\begin{multline}\label{eq:JOsFEMequi}
    \J(\Omega_s) %
     = \sum_{i \in \Ci}\int_{\Kref}  \left(u_{(\F + s\VV \circ \F)\circ\VF_i^{-1}} \circ (\F + s\VV \circ \F)\right) \\
      \vert\det(\VD (\F + s\VV \circ \F))\vert\DXref\,.
\end{multline}
{Let $\{g_i\}_{i\in\Ci}$ be the collection of maps defined by
\begin{equation*}
g_i(\VT) \coloneqq (u_{\VT\circ\VF_i^{-1}}\circ\VT)\vert\det(\VD\VT)\vert\,.
\end{equation*}
Then, formula \eqref{eq:JOsFEMequi} can be rewritten as
\begin{equation*}
\J(\Omega_s) =  \sum_{i \in \Ci}\int_{\Kref}  g_i(\F + s\VV \circ \F)\DXref\,,
\end{equation*}
and}
taking the derivative of \eqref{eq:JOsFEMequi} with respect to $s$ implies
\begin{multline} \label{eq:dJgen-ufl}
\dJ(\Omega)[\VV] = \sum_{i \in \Ci}\int_{\Kref} \ders{g_i(\F + s\VV \circ \F)}\DXref\,.
\end{multline}
Equation \eqref{eq:dJgen-ufl} gives an alternative {and equivalent} expression
for the shape derivative \eqref{eq:dJgen}.
{However, to derive formula \eqref{eq:dJgen} it is necessary
to follow shape calculus rules by hand, which is often a tedious and error prone
exercise. Equation \eqref{eq:dJgen-ufl}, by contrast,
can be derived automatically with finite element software.
Indeed,
to evaluate $\J(\Omega)$, standard finite element software rewrites it as
\begin{equation*}
\J(\Omega) = \sum_{i \in \Ci}\int_{\Kref} g_i(\F)\,.
\end{equation*}
{In UFL the maps $\{g_i\}_{i\in \Ci}$ are constructed symbolically and in an automated fashion.}
Therefore, it is possible to evaluate $\dJ(\Omega)[\VV]$ by performing the
steps necessary for the assembly of $\J(\Omega)$ and, at the appropriate time,
differentiating the maps $\{g_i\}_{i\in\Ci}$. To be precise, this differentiation
corresponds to a standard G\^{a}teaux directional derivative, because
the integrand in \eqref{eq:dJgen-ufl} corresponds to the following limit
\begin{equation*}
\ders{g_i(\F + s\VV \circ \F)}=\lim_{s\searrow 0}\frac{g_i(\F + s\VV \circ \F)-g_i(\F)}{s}\,,
\end{equation*}
which can be interpreted as the G\^{a}teaux directional derivative of $g_i$
at  $\VT=\F$ in the direction $\VV\circ\F$ \cite[Def. 1.29]{HiPiUlUl09}.
This viewpoint is important to correctly implement this differentiation step in the existing pipeline in UFL
(see Figure \ref{fig:UFLworkflow}).}
We emphasize that this also enables computing higher-order shape
derivatives by simply taking higher-order G\^{a}teaux derivatives
in \eqref{eq:dJgen-ufl}.

\begin{remark}
  Lagrange finite element global basis functions are obtained by
  gluing local parametric basis functions, that is, basis functions
  $\{b_m^i\}_{m\in\Cm}$ defined only on $K_i$ and of the form
  $b_m = \myref{b_m}\circ \F^{-1}$, where $\{\myref{b_m}\}_{m\in\Cm}$
  is the set of reference local basis functions, which are defined
  only on the reference element $\Kref$.  If $\VV$ lives in a Lagrange
  finite element space built on the partitioning $\{K_i\}_{i\in \Ci}$,
  it is possible to evaluate $\dJ(\T)[\VV]$ by computing the
  G\^{a}teaux derivative in \eqref{eq:dJgen-ufl} in the direction of
  the reference local basis functions $\{\myref{b_m}\}_{m\in\Cm}$
  (instead of in the direction $\VV\circ\F$) and summing these
  values. This allows us to fully rely on the symbolic
  differentiation capabilities of UFL.
\end{remark}
\begin{remark}
  The approach does not rely on the element being affinely mapped, but extends to elements that are mapped using a Piola transform such as the Raviart-Thomas or Nedelec elements.
  However, it does not work for elements such as the Hermite element that require different pullbacks for point evaluation and derivative degrees of freedom.
\end{remark}

\section{Examples}
\label{sec:examples}
In this section, 
we consider three examples based on \eqref{eq:proto-j} that cover most applications.
For these examples, we give explicit expressions of $\dJ$ using
\eqref{eq:dJgen} and \eqref{eq:dJgen-ufl}
and show how to compute $\dJ$
with the finite element software
Firedrake\footnote{Examples using FEniCS \cite{Logg:2010,Logg:2012} will be almost identical,
modulo small differences in setting up initial conditions} \cite{Rathgeber:2016}.
To shorten the notation, we define $\VV_i \coloneqq\VV\circ\VF_i$.
\begin{example}
  \label{example:1}
  Let the integrand be independent of $\Omega$, i.e. $u_{\VP_s}= u$
  for some function $u$.  Then, the chain rule implies that
  $\ders{u_{\VP_s}\circ\VP_s}=\ders{u\circ\VP_s}=\nabla u\cdot \VV$.
  Recalling $\ders{\det\VD\VP_s} = \Div(\VV)$, we conclude that
  equation \eqref{eq:dJgen} becomes
  \begin{equation}
    \label{eq:dJindep} \dJ(\Omega)[\VV] = \int_{\Omega}  \nabla u \cdot \VV +  u\Div(\VV)\DX\,.  
  \end{equation}
  On the other hand, inserting $u_{\VP_s}=u$ into
  \eqref{eq:dJgen-ufl}, we obtain the equivalent
  expression: %
  \begin{multline}\label{eq:dJindep-ufl}
    \dJ(\Omega)[\VV] %
    = \sum_{i \in \Ci}\int_{\Kref}  
    \big(\VV_i \cdot (\nabla u \circ \F) \\+ (u\circ \F) \Tr(\VD \VV_i\VD \F^{-1})\big)\vert\det(\VD\F)\vert \DXref\,.
  \end{multline}
  Example code is shown in \textcolor{black}{Listing} \ref{lst:example1}.
  {
      Functionals with domain independent integrands are used in
      applications including image segmentation \cite{hintermuller_second_2004} or, when $u\equiv 1$, to enforce volume constraints in shape optimization.}

  \begin{mylisting}[htb!]
  \hspace{-1.52em}
\includegraphics[width=1.06\linewidth]{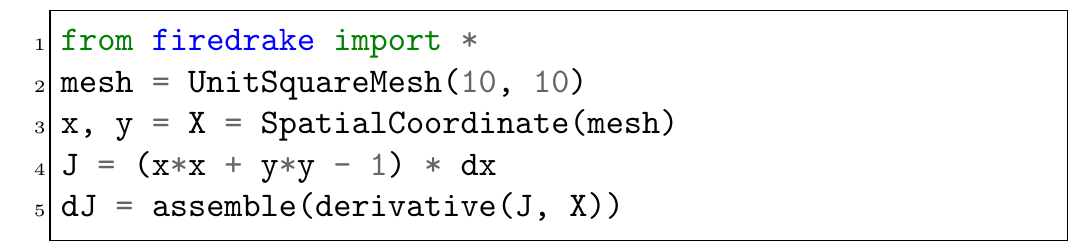}
  \caption{%
  \textcolor{black}{Firedrake code to compute $\dJ$ from Example \ref{example:1} when {$u(x,y) = x^2 + y^2 - 1$.}}}
  \label{lst:example1}
  \end{mylisting}

\end{example}

\begin{example}
  \label{example:2}
  Let $\{V_h(\Omega_s)\}_s$ be a family of scalar finite element
  spaces such that the global basis functions $\{B_s^i\}_i$ of
  $V_h(\Omega_s)$ are of the form $B_s^i = B^i\circ\VP_s^{-1}$, where
  $\{B^i\}_i$ are basis functions of $V_h(\Omega)$ {and $\VP_s^{-1}$
  is the inverse of $\VP_s$, that is, $\VP_s^{-1}(\VP_s(\Vx))= \VP_s(\VP_s^{-1}(\Vx))=\Vx$
  for every $\Vx\in \bbR^d$}.  Let
  $v_{\VP_s}\in V_h(\Omega_s)$ and
  $u_{\VP_s}=v_{\VP_s}+\| \nabla v_{\VP_s} \|^2$.  Since
  \begin{equation}
    v_{\VP_s} \circ \VP_s + \|(\nabla v_{\VP_s})\circ \VP_s\|^2 = v_{\VP_0} + \|\VD\VP_s^{-T} \nabla v_{\VP_0}\|^2\,,
  \end{equation}
  \textcolor{black}{equation} \eqref{eq:dJgen} becomes
  \begin{multline}
    \dJ(\Omega)[\VV]
    = \int_\Omega (v_{\VP_0} + \| \nabla v_{\VP_0}\|^2) \Div(\VV)
   \\ -2 \nabla v_{\VP_0} \cdot (\VD\VV^{T} \nabla v_{\VP_0})  \DX\,.
  \end{multline}
  On the other hand, note that for any $\hat\Vx\in \hat K$ and for
  $i\in\Ci$, it holds $v_{\VP_s}(\VP_s(\F(\xref))) = v_i(\xref)$,
  where $v_i$ is a linear combination of the local basis functions
  $\{\myref{b_m}\}_{m\in\Cm}$ defined on the reference element
  $\hat K$.  Therefore,
  \begin{multline}
    v_{\VP_s} \circ \VP_s \circ \F + \|(\nabla v_{\VP_s})\circ \VP_s \circ \F\|^2 \\
    = \hat v_i + \|\VD(\VP_s \circ \F)^{-T} \nabla \hat v_i\|^2
    \quad \text{on } \hat K\,,
  \end{multline}
  and \eqref{eq:dJgen-ufl} becomes
  \begin{equation}
    \begin{aligned}
      \dJ(\Omega)[\VV] 
      = \sum_{i \in \Ci} \int_{\Kref} \big((\hat v_i + \|\VD\F^{-T}
        \nabla v_i\|^2) \Tr(\VD\VV_i \VD\F^{-1}) \\
        - 2 \nabla \hat v_i
        \cdot \VD\F^{-1}\VD\VV_i\VD\F^{-1}\VD\F^{-T}\nabla \hat
        v_i\big) \vert\det(\VD\F)\vert \DXref.
    \end{aligned}
  \end{equation}
  \textcolor{black}{Listing} \ref{lst:example2} shows code for this case, using as $v_{\VP_0}$
  the piecewise affine Lagrange interpolant of $\sin(x)\cos(y)$.
\begin{mylisting}[htb!]
\hspace{-1.52em} \includegraphics[width=1.06\linewidth]{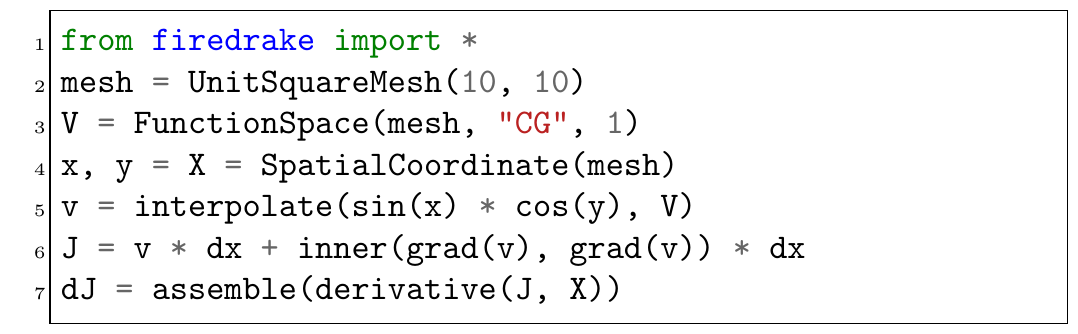}
\caption{%
\textcolor{black}{Firedrake code to compute $\dJ$ from Example \ref{example:2}.}
Note that in this
case, $v$ does not depend explicitly on $x$ and $y$.}
\label{lst:example2}
\end{mylisting}
\end{example}

\begin{example}
  \label{example:3}
  Let $u_{\VP_s}$ be the finite element solution to the boundary value
  problem
  \begin{equation}\label{eq:PDEconstr}
    -\Delta u_{\VP_s} + u_{\VP_s} = f \quad \text{in } \Omega_s\,,
    \quad\nabla u_{\VP_s} \cdot\Vn = 0\quad \text{on } \partial \Omega_s\,.
  \end{equation}
  In this case, the functional \eqref{eq:proto-j} is said to be
  PDE-constrained, and computing its shape derivative is less
  straightforward. The standard procedure is to introduce an
  appropriate Lagrangian functional \cite[Ch.~10, Sect.~5]{DeZo11}.
  For this example, the Lagrangian is
  \begin{equation}\label{eq:Lagrangian}
    \mathrm{L}_s(u_{\VP_s},v_s) \coloneqq \J(\Omega_s) + \mathrm{e}_s(u_{\VP_s},v_s),
  \end{equation}
  where
  \begin{equation}\label{eq:state}
    \mathrm{e}_s(u_{\VP_s}, v_s) \coloneqq 
    \int_{\Omega_s}\nabla u_{\VP_s}\cdot \nabla v_s + u_{\VP_s}v_s - fv_s\DX\,%
  \end{equation}
  stems from the weak formulation of the PDE constraint
  \eqref{eq:PDEconstr}.  The shape derivative $\dJ$ is equal to the
  shape derivative of
  $\mathrm{L}_s(u\circ\VP_s^{-1},p\circ\VP_s^{-1})$, where $u$ is the
  solution to \eqref{eq:state} for $s=0$ and $p\in V_h(\Omega)$ is the
  solution to an adjoint boundary value problem.  The shape derivative
  of $\mathrm{L}_s(u\circ\VP_s^{-1},p\circ\VP_s^{-1})$ can be computed
  as in \textcolor{black}{Example} \ref{example:2}.  The result is
  \begin{multline}\label{eq:dJconstr}
    \dJ(\Omega)[\VV] = \int_{\Omega} (u + \nabla u\cdot \nabla p + up - fp)\Div(\VV)\\
    - p\nabla f\cdot\VV
    -  \nabla u\cdot (\VD\VV+\VD\VV^{T})\nabla p\DX\,.
  \end{multline}
  For this example, we omit the equivalent formula on the reference
  element because of its length.  However, as \textcolor{black}{Listing} \ref{lst:example3}
  shows, UFL removes the tedium of deriving the shape derivative, and
  we can easily compute $\dJ$.
\begin{mylisting}[htb!]
\hspace{-1.52em} \includegraphics[width=1.06\linewidth]{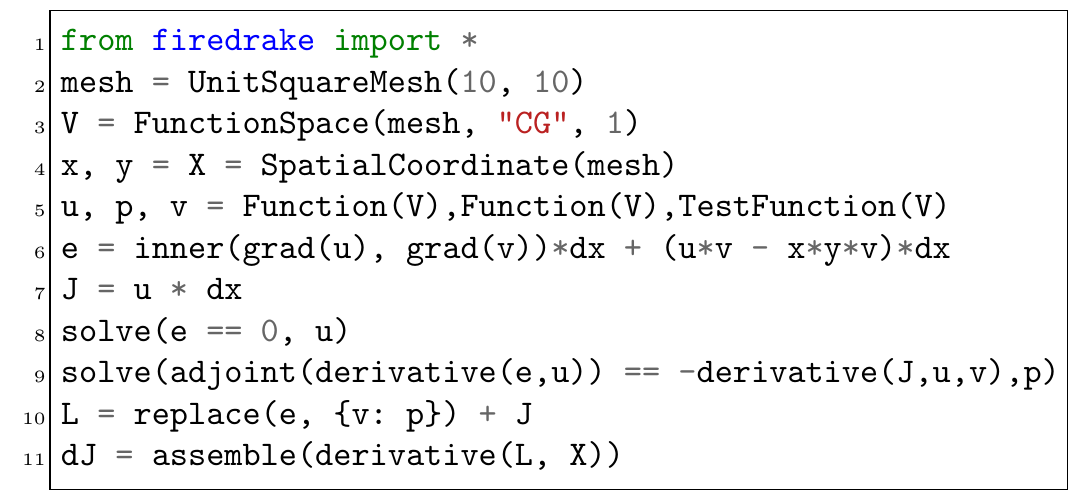}
\caption{%
\textcolor{black}{Firedrake code to compute $\dJ$ from Example \ref{example:3} when $f(x,y) = xy$ in \eqref{eq:state}.}}
\label{lst:example3}
\end{mylisting}
\end{example}
\begin{remark}
  With appropriate modifications, the same code can be use for
  functionals constrained to boundary value problems with Neumann or
  Dirichlet boundary conditions.  For the Neumann case, it is
  sufficient to add the Neumann forcing term in line 6 of
  \textcolor{black}{Listing }\ref{lst:example3}. For the Dirichlet case, one needs to replace
  \texttt{u} with \texttt{u+g} in lines 6 and 7 (where \texttt{g} is
  the function defined in terms of \texttt{X} that describes the
  Dirichlet boundary condition) and impose homogeneous Dirichlet
  boundary conditions in lines 8 and 9.
\end{remark}
\begin{remark}
  To evaluate the action of the shape Hessian of a PDE-constrained
  functional, one can follow the instructions given in
  \cite[p. 65]{HiPiUlUl09}.  Note that by computing shape derivatives
  as in \eqref{eq:dJgen-ufl}, it is straightforward to combine shape
  derivatives of $\mathrm{L}_s(u\circ\VP_s^{-1},p\circ\VP_s^{-1})$
  with standard G\^{a}teaux derivatives with respect to
  $u\circ\VP_s^{-1}$.
\end{remark}

\section{Code validation}
\label{sec:implementation-validation}
\begin{figure}[htb!]\label{fig:Taylor}
    \begin{center}
        \includegraphics{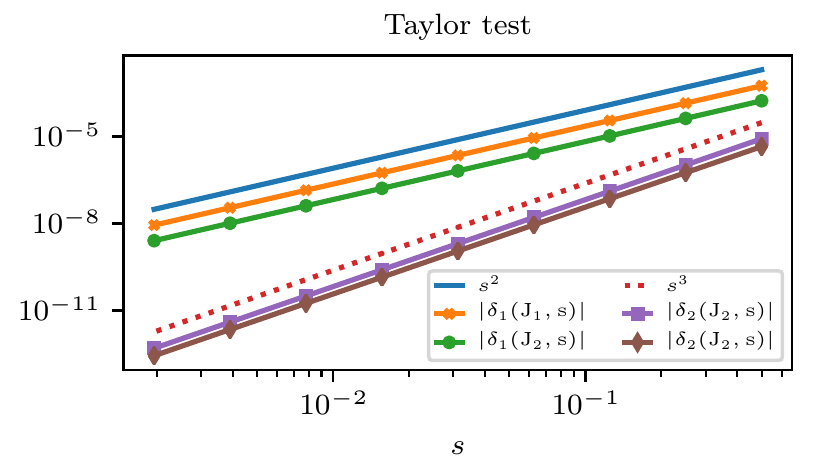}
        \caption{Taylor test for Examples 1 and 3. The convergence rates match the expected convergence.}
    \end{center}   
\end{figure}

We validate our implementation by testing
that the Taylor expansions truncated to first and second
order satisfy the asymptotic conditions
\begin{equation}\label{eq:Taylor}
\rm\delta_1(\J, s) = O(s^2) \quad \text{and} \quad \rm\delta_2(\J, s) = O(s^3)\,,
\end{equation}
where
\begin{equation*}
\rm\delta_1(\J, s) \coloneqq \|\J({\Omega_s}) - \J(\Omega) - s \dJ(\Omega)[\VV]\|
\end{equation*}
and
\begin{multline*}
\rm\delta_2(\J, s) \coloneqq \|\J({\Omega_s}) - \J(\Omega)\\
- s \dJ(\Omega)[\VV] - \frac{1}{2} s^2 \ddJ(\Omega)[\VV, \VV]\|\,.
\end{multline*}
In \textcolor{black}{Figure }\ref{fig:Taylor}, we plot the values of $\rm\delta_1$ and $\rm\delta_2$
for $s=2^{-1}, 2^{-2},$ $\dots, 2^{-10}$, and $\J$
as in \textcolor{black}{Examples 1 and 3} %
from the previous section (we
denote these functionals $\J_1$ and $\J_2$ respectively).
The vector field $\VV$ is chosen randomly. This experiment clearly
displays the asymptotic rates predicted by \eqref{eq:Taylor}.
We have repeated this numerical experiment
for many other test cases, including functionals that are not linear
in $u$, {functionals given by integrals over $\partial\Omega$ involving the normal $\Vn$,}
and functionals that are constrained to linear and nonlinear
boundary value problems with nonconstant
right-hand sides and nonconstant Neumann and Dirichlet boundary conditions.
In every instance, we
have observed the asymptotic rates predicted by \eqref{eq:Taylor}.
The code for these numerical experiments is available at \cite{zenodoFiredrake20190402.0}.

\section{Shape optimization of a pipe}
\label{sec:shapeoptexample}
In this section, we show how to use Firedrake and the new UFL capability to code
a PDE-constrained shape optimization algorithm. As test case, we consider the
optimization of a pipe to minimize the dissipation of kinetic energy of the fluid into heat.
This example is taken from \cite[Sect. 6.2.3]{SchThesis}. To simplify the exposition,
we use a very simple optimization strategy. At the end of the section, we will comment
on possible improvements.

\begin{figure}
  \includegraphics[width=\linewidth]{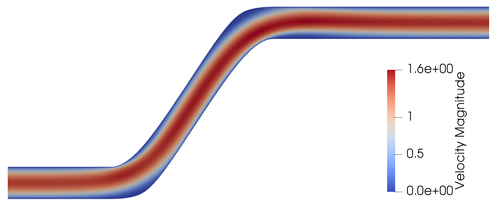}
  \includegraphics[width=\linewidth]{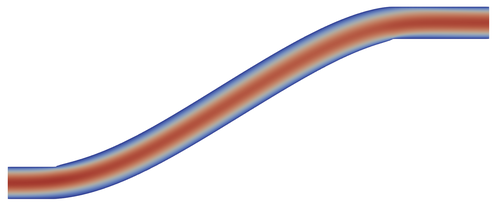}
  \caption{Initial (top) and optimized (bottom) shape of a pipe connecting a given inflow and outflow.}
  \label{fig:pipe}
\end{figure}

The initial design of the pipe is shown in Figure \ref{fig:pipe} (top).
The pipe contains viscous fluid (with viscosity $\mu$), which flows in from the left and is
modelled using the incompressible Navier-Stokes equations.
To be precise, let $\Omega$ be the shape of the pipe, $\Gamma\subset \partial\Omega$
be the outflow boundary of the pipe (that is, the end of the pipe on the right),
and $\Vu$ and $p$ be the velocity
and the pressure of the fluid, respectively. Then, $\Vu$ and $p$ satisfy
\begin{align*}
-\nu \Delta \Vu + \Vu \nabla \Vu + \nabla p &= 0 & \text{in } \Omega\,,\\
\Div \Vu &= 0& \text{in } \Omega\,,\\
\Vu &= \Vg &\text{on } \partial\Omega\setminus \Gamma\,,\\
p\Vn - \nu \nabla u\cdot \Vn & = 0 &\text{on } \Gamma\,.
\end{align*}
Here $\Vg$ is given by a Poiseuille flow at the inlet and is zero on the walls of the pipe

\begin{figure}[htb!]
  \includegraphics{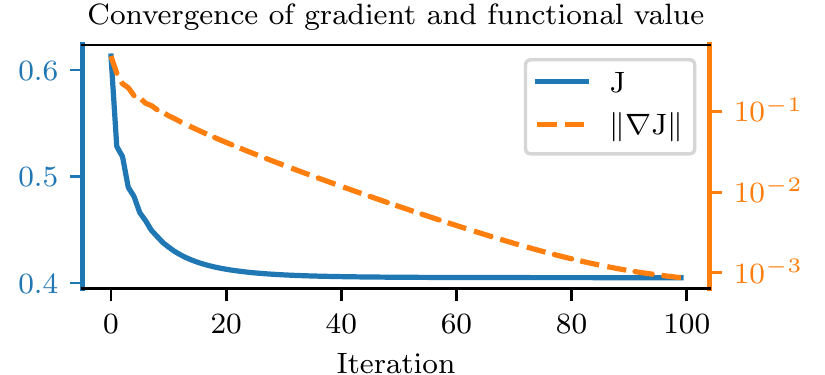}
  \caption{The value of the objective (plotted in linear scale) is reduced from approximately $0.61298$ to $0.40506$.
  The $H^1$-norm of the gradient (plotted in logarithmic scale) is reduced from $0.487274$ to $0.000870$.}
  \label{fig:pipe-convergence}
\end{figure}

The goal is to modify the central region of the pipe
so that the shape functional
\begin{equation*}
\J(\Omega) = \int_\Omega \nu \nabla\Vu : \nabla\Vu \DX
\end{equation*}
is minimized.
To solve this shape optimization problem,
we parametrize the initial design with a polygonal mesh and update
the node coordinates using a descent direction optimization algorithm
with fixed step size.
As descent directions, we use Riesz representatives of the shape gradient
with respect to the inner product induced by the Laplacian{, i.e. at each step the deformation is given by the solution to
\begin{equation}
    \begin{aligned}
        -\Delta \VV &= -\dJ(\Omega) && \text{in } \Omega\\
        \VV &= 0 && \text{on fixed boundaries.}
    \end{aligned}
\end{equation}
This approach is also known as Laplace smoothing.
}
To avoid degenerate results, we penalize changes
of the pipe volume.
The whole algorithm, comprising of state and adjoint equations and
shape derivatives, is contained in Listing \ref{lst:pipe} and described in
detail in the following paragraph. The optimized shape is displayed in Figure \ref{fig:pipe} (bottom),
the convergence history is in Figure \ref{fig:pipe-convergence}.
These results are compatible with those in \cite[Sect. 6.2.3]{SchThesis} and
clearly indicate the success of the shape optimization algorithm.

\paragraph{Description of Listing \ref{lst:pipe}}
In line 2-4, we load the finite element mesh \texttt{pipe.msh} and extract the vertex coordinates.
This mesh is generated with Gmsh \cite{geuzaine2009gmsh} and is available as part of \cite{zenodoFiredrake20190402.0}.
Lines 5-8 define the Gramian matrix of the inner product employed
to compute descent directions. In lines 9-14, we define the space of P2-P1
Taylor-Hood finite elements, which we use to discretize the weak formulation of
the Navier-Stokes equations, and set up the functions containing the solutions
to the state and adjoint equation as well as the test functions for the weak form.
In lines 15-22 we define the weak formulation of the Navier-Stokes equations and
certain parameters to prescribe the use of the MUMPS direct solver \cite{amestoy2000mumps} to solve
the linearized equations.  In lines 23-29, we define the shape functional
$\J$, the functional describing the volume of the shape, as well as the
Lagrangian and its derivative. {In particular, note that
the shape derivative of the Lagrangian can be computed with
the simple command \verb!dL = derivative(L, X)! in line 28.
Without the new automatic shape differentiation capability in UFL,
line 28 would have to be replaced with the following formula
\begin{verbatim}
 dL = -inner(nu*grad(u)*grad(W), grad(v))*dx
      -inner(nu*grad(u), grad(v)*grad(W))*dx
      -inner(v, grad(u)*grad(W)*u)*dx
      +tr(grad(v)*grad(W))*p*dx
      -tr(grad(u)*grad(W))*q*dx
      +div(W)*inner(nu*grad(u), grad(v))*dx
      -div(W)*inner(div(v), p)*dx
      +div(W)*inner(div(u), q)*dx
      +div(W)*inner(v, grad(u)*u)*dx
      +nu*inner(grad(u), grad(u))*div(W)*dx
      -2*nu*inner(grad(u)*grad(W), grad(u))*dx
\end{verbatim}}
\noindent In lines 30-35 we set up a function that updates the solution to the state and
the adjoint equations. 
{We emphasize that this shape optimization problem is not self-adjoint and that UFL derives the adjoint equation automatically.}
Note that, whenever the function
\texttt{solve\_state\_and\_adjoint} is called, the new values of the velocity $\Vu$
are stored in the file \texttt{u.pvd} (which can visualized using Paraview \cite{ahrens2005paraview}).
Finally, lines 36-46 contain the optimization algorithm: for 100 iterations we
compute the shape derivative and penalize volume changes (lines 38-40), compute the descent direction (lines 41-42),
update the domain (line 45), and
update the state and adjoint solutions (line 46). 

\begin{mylisting}[htb!]
\hspace{-1.52em} \includegraphics[width=1.06\linewidth]{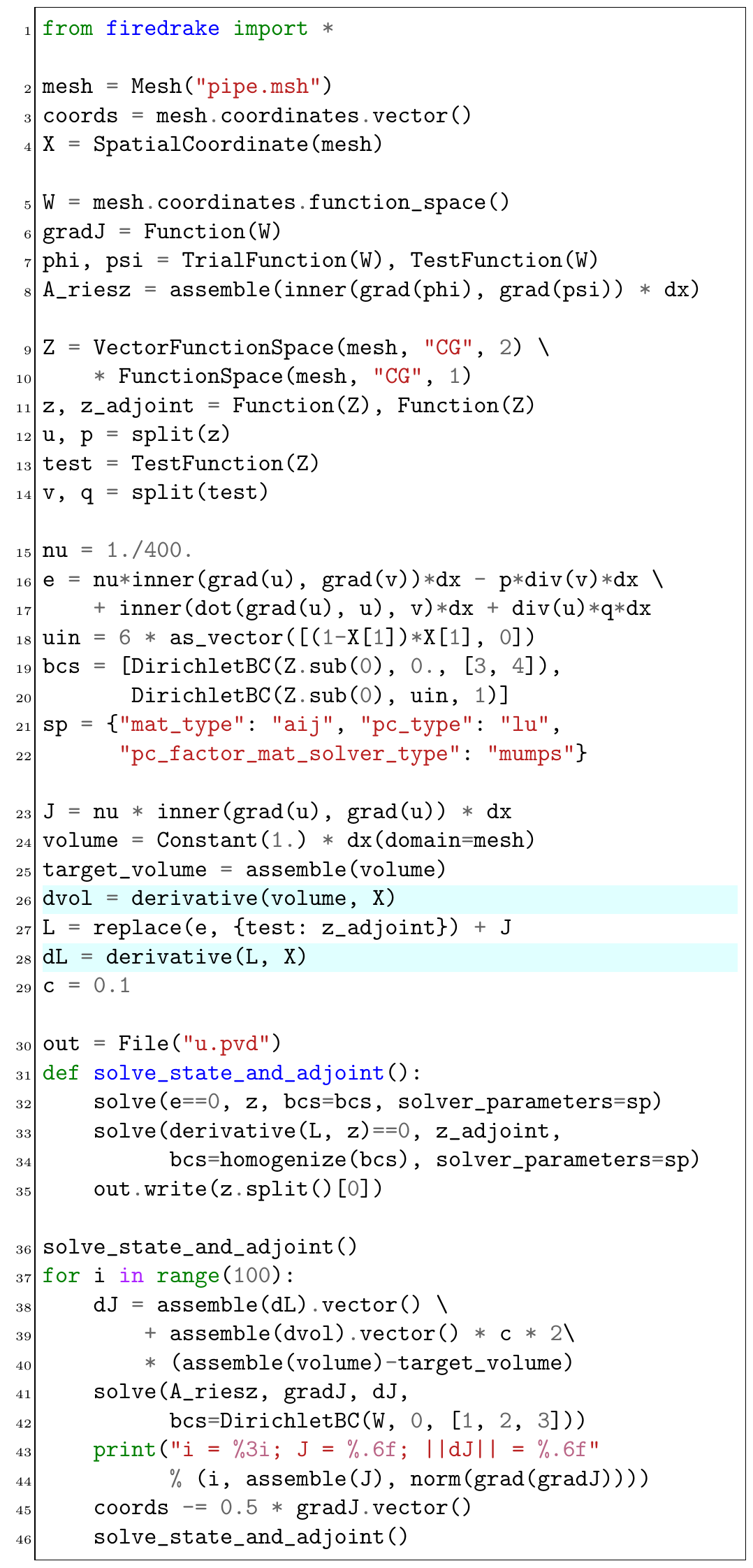}
\caption{Firedrake code to optimize the shape of a pipe and minimize the dissipation of kinetic energy into heat. \textcolor{black}{Lines 26 and 28 use the newly developed automatic shape differentiation.}}
\label{lst:pipe}
\end{mylisting}

\begin{remark}
The optimization algorithm of Listing \ref{lst:pipe} is based on a simple optimization
strategy and can be improved in several ways, at the mere cost of adding lines of code.
For instance, instead of using a fixed step-size and a fixed number of iterations,
one could implement an adaptive step-size selection and stopping criteria.
Additionally, one could experiment with different inner products to define descent
directions \cite{IgStWe17}, as well as compute second order derivatives of
$\J$ and implement (quasi-)Newton methods \cite{Sc18}. Despite the room for improvement,
we would like to stress that Listing \ref{lst:pipe} can be readily used for a 3D problem by simply passing a 3D mesh
and changing the inflow boundary condition in line 18.
\end{remark}

\section{Discussion}
\label{sec:discussion}
{We have presented a new and equivalent formulation of shape derivatives in the context of
finite elements as G\^ateaux-derivatives on the reference element.  While the
formulation applies to finite elements in general, we have implemented this
new approach in UFL due to its extensive support for symbolic calculations.
This new UFL capability allows computing shape derivatives of
functionals that are defined as volume or boundary integrals, and that are constrained to linear and non-linear PDEs.}
During shape differentiation, our code treats finite element functions
and global functions differently. This behavior is correct and
necessary to handle PDE-constraints properly.
In combination with a finite element software package, such as FEniCS
or Firedrake, that takes as input UFL, this enables the
entirely automated shape differentiation of functionals subject to boundary value problems.
\textcolor{black}{This notably simplifies tackling PDE-constrained shape optimization problems.}

Compared to the existing shape differentiation toolbox FEMorph,
our code does not compute shape derivatives in strong form
because it neither relies on shape calculus differentiation rules
nor performs integration by parts.
However, in practice we do not consider this a limitation as it has been shown in
\cite{HiPaSa15,Be10,Zh18} that the weak form is superior when the state and
the adjoint equations are discretized by finite elements.

\begin{acknowledgements}
DAH is supported by the Natural Environment Research Council [grant number
NE/K008951/1]. LM is supported by the Engineering and Physical Sciences
Research Council [grant number EP/L000407/1].  FW is supported by the EPSRC
Centre For Doctoral Training in Industrially Focused Mathematical Modelling
[grant number EP/L015803/1].\\~\\
\textbf{Contributions}
This work originated in a discussion between the four authors at the FEniCS 18 conference, where AP and DAH suggested to calculate shape derivates as in Section \ref{sec:theory}.
FW implemented this idea in UFL with help from LM.
The manuscript was written by AP and FW, with feedback from LM and DAH.
\\~\\
\textbf{Conflict of interest statement}
 On behalf of all authors, the corresponding author states that there is no conflict of interest.
 \\~\\
\textbf{Replication of results}
The code for the numerical experiments is available at \cite{zenodoFiredrake20190402.0}.
\end{acknowledgements}

\bibliographystyle{plainurl}

\end{document}